\documentclass[11pt]{article}
\usepackage{anysize}\marginsize{3.5cm}{3.5cm}{1.3cm}{2cm}
\usepackage{amssymb}
\newtheorem{satz}{Satz}[section]

\newtheorem{corollary}[satz]{Corollary}

\newtheorem{example}[satz]{Example}

\newtheorem{proposition}[satz]{Proposition}
\newtheorem{remark}[satz]{Remark}

\newtheorem{theorem}[satz]{Theorem}

\newcommand\beginproof[1]{\trivlist\item[\hskip\labelsep{\em #1.}]}

\newcommand\proofof[1]{\beginproof{Proof of #1}}
\def\endproof{\hspace*{\fill}\endproofsymbol\endtrivlist}

\def\endproofsymbol{\frame{\rule[0pt]{0pt}{6pt}\rule[0pt]{6pt}{0pt}}}

\renewcommand\ge{\geqslant}  
\renewcommand\le{\leqslant}  
\renewcommand\geq{\geqslant}  
\renewcommand\leq{\leqslant}  
\renewcommand\epsilon{\varepsilon}
\renewcommand\phi{\varphi}

\renewcommand\O{{\cal O}}
\renewcommand\P{\mathbb P}

\newcommand\be{\begin{eqnarray*}}
\newcommand\ee{\end{eqnarray*}}

\newcommand\eqnref[1]{(\ref{#1})}
\newcommand\eqdef{\stackrel{\mbox{\scriptsize def}}{=}}
\newcommand\eps{\varepsilon}
\newcommand\tensor{\otimes}

\newcommand\rounddown[1]{\left\lfloor#1\right\rfloor}
\newcommand\tfrac[2]{{\textstyle\frac{#1}{#2}}}
\newcommand\isom{\cong} 
\newcommand\inparen[1]{\textnormal{(}{#1}\textnormal{)}}

\newcommand\mult{{\rm mult}}

\newcommand\newop[2]{\newcommand#1{\mathop{\rm #2}\nolimits}}
\newop\Sym{Sym}

\begin{document}

\title{Seshadri constants and the generation of jets}
\author{Thomas Bauer and Tomasz Szemberg}
\date{February 6, 2009}
\maketitle
\thispagestyle{empty}

\section*{Introduction}
\begin{plaintheoremnumbers}

   Consider a smooth projective variety $X$ and an ample line
   bundle $L$ on $X$.
   The \textit{Seshadri constant}
   of $L$
   at a given point $x\in X$ is the real number
   $$
      \eps(L,x)\eqdef\inf_C\frac{L\cdot C}{\mult_xC} \,,
   $$
   where the infimum is taken over all irreducible
   curves $C$ passing through $x$.
   It is well-known that
   $\eps(L,x)$ encodes
   asymptotic information about
   the linear series $|kL|$ for $k\gg 0$.
   Specifically, denote for $k\ge 1$
   by $s(kL,x)$ the maximal integer $s$
   such that the linear series $|kL|$ generates $s$-jets at $x$,
   i.e., the maximal integer $s$ such that the evaluation map
   $$
      H^0(X,kL)\longrightarrow H^0(X,kL\tensor\O_X/\mathcal I_x^{s+1})
   $$
   is onto.
   Then one has
   $$
      \eps(L,x)=\lim_{k\to\infty} \frac{s(kL,x)}k
   $$
   (see \cite[6.3]{Dem92}).
   So if the numbers $s(kL,x)$ were known for $k\gg 0$, then one
   could compute the Seshadri constant as a limit of a
   sequence of rational numbers.
   In all situations that we know of, however, Seshadri
   constants have not been determined in this way, but rather
   by finding
   suitably singular curves.
   It is therefore interesting to ask:
   \begin{quote}\it
      When the value $\eps(L,x)$ is known, what can be said about the
      numbers $s(kL,x)$ for $k\ge 1$ ?
   \end{quote}
   It seems to us that this question is close in spirit to
   Demailly's original purpose when defining Seshadri constants
   in \cite{Dem92}.

   Our first result shows
   that under certain non-positivity assumptions on the
   canonical divisor of $X$ the range for the numbers $s(kL,x)$ is
   quite restrictive.

\begin{theorem}\label{range}
   Let $X$ be a smooth projective variety of dimension $n$.

   \begin{itemize}
      \item[\rm (a)] If $X$ is Fano, then for any integer $k\ge 1$
      $$
      \rounddown{(k+1)\,\eps(-K_X,x)}-(n+1)\le
      s(k(-K_X),x)\le\rounddown{k\,\eps(-K_X,x)}\,.
      $$
   \item[\rm (b)] If $K_X=\O_X$, then for any ample line bundle $L$ on
      $X$, any point $x\in X$, and any integer $k\ge 1$ we have
      $$
      \rounddown{k\,\eps(L,x)}-(n+1)\le
      s(kL,x)\le\rounddown{k\,\eps(L,x)}\,.
      $$
   \end{itemize}
\end{theorem}

   Note that the upper bound -- which in fact holds in
   both parts
   without any assumptions on $X$ -- is well-known and is
   stated here merely for the sake of completeness.
   We will show in Proposition~\ref{cones}
   that $s(kL,x)$ can attain this upper bound
   only if $\eps(L,x)$ is either computed by a
   smooth curve or if it is not computed by a curve at all.
   As for the lower bounds, we will prove
   somewhat stronger statements in Propositions
   \ref{prop-fano} and \ref{prop-K-zero} respectively.

   The above bounds on the numbers $s(kL,x)$ lead in particular to the
   following new characterization of projective spaces.

\begin{theorem}\label{char-pn}
   Let $X$ be a smooth Fano variety of dimension $n$ such that there
   exists a point $x\in X$ with
   $$\eps(-K_X,x)=n+1.$$
   Then $X$ is the projective space $\P^n$.
\end{theorem}

   We will in fact show that on Fano varieties different from
   $\P^n$ one has $\eps(-K_X,x)\le n$ for all points $x$
   (Theorem~\ref{better_than_pn}).

\end{plaintheoremnumbers}

\section{Fano varieties}

   We begin by considering jets of the anticanonical bundle on Fano
   varieties.

\begin{proposition}\label{prop-fano}
   Let $X$ be a smooth Fano variety of dimension $n$, and let
   $x\in X$. Let $\eps=\eps(-K_X,x)$.

   If $\eps<\sqrt[n]{(-K_X)^n}$
   or if $\sqrt[n]{(-K_X)^n}$ is not an integer, then
   $$
      \rounddown{(k+1)\,\eps}-n\le
      s(k(-K_X),x)\le\rounddown{k\,\eps} \,.
   $$
   In the alternative case, one has
   $$
      (k+1)\eps-(n+1)\le s(k(-K_X),x)\le k\eps \,.
   $$
\end{proposition}
\begin{proof}
   The upper bound -- which in fact holds without any assumption
   on $K_X$ -- follows
   from the fact that $\eps(L,x)$ is not only the limit, but
   also the
   \textit{supremum} of the
   numbers $\frac1k s(kL,x)$, see \cite[(6.3)]{Dem92}.

   The lower bound is proven via vanishing as in
   \cite[Proposition 5.1.19(i)]{PAG}: One considers
   the blow-up $f:X'\to X$ of $X$ at $x$ with exceptional divisor
   $E$ over $x$. For $|k(-K_X)|$ to generate $s$-jets at $x$
   it is enough to have
   $H^1(kf^*(-K_X)-(s+1)E)=0$.
   As $K_{X'}=f^*(K_X)+(n-1)E$, this vanishing will follow if the line
   bundle $(k+1)f^*(-K_X)-(s+n)E$ is nef and big.
   Let us write $\eps=\eps(-K_X,x)$, and
   suppose first that $\eps<\sqrt[n]{(-K_X)^n}$. In that case
   $(k+1)f^*(-K_X)-(s+n)E$ is nef and big
   as long as
   $\frac{s+n}{k+1}\le\eps$.
   Therefore
   \begin{equation}\label{goodlower}
      s(k(-K_X),x)\ge\rounddown{(k+1)\eps}-n \ .
   \end{equation}
   Suppose then that $\eps=\sqrt[n]{(-K_X)^n}$. In that case
   the line bundle in question is ample
   if $\frac{s+n}{k+1}<\eps$, i.e, if $s<(k+1)\eps-n$.
   Now, if $\eps$ is not an integer, then this inequality is
   equivalent to $s\le\rounddown{(k+1)\eps}-n$, so that we
   get again
   \eqnref{goodlower}. Finally, if $\eps=\sqrt[n]{(-K_X)^n}$
   is an integer, then we get
   $$
      s(k(-K_X),x)\ge\rounddown{(k+1)\eps}-(n+1)=(k+1)\eps-(n+1) \,,
   $$
   as claimed.
\end{proof}

\begin{example}[Projective Space]\rm
   For $X=P^n$ and $L=\O_{\P^n}(1)$ one has $\eps(L,x)=1$ for all
   $x\in X$ and $s(kL,x)=k$ for all $k$. So here the value
   $$s(k(-K_X),x)=k(n+1)$$
   lies at the upper bound given by Proposition~\ref{prop-fano}. We
   will show that
   projective spaces are the only Fano varieties where this
   happens (Theorem~\ref{better_than_pn}).
\end{example}

   The proposition leads immediately to a surprising upper bound
   on Seshadri constants:

\begin{corollary}\label{cor-fano}
   For any smooth Fano variety $X$, one has
   $$
      \eps(-K_X,x)\le n+1
   $$
   for all $x\in X$.
   If $\eps(-K_X,x)<\sqrt[n]{(-K_X)^n}$
   or if $\sqrt[n]{(-K_X)^n}$ is not an integer,
   then the stronger
   inequality
   $$
      \eps(-K_X,x)\le n
   $$
   holds.
\end{corollary}

\begin{proofof}{the corollary}
   This follows from the inequalities in the proposition. In
   fact,
   it is enough to show that if $\eps$ is a positive real number
   and
   $b$ is an integer such that
   \begin{equation}\label{round-ineq}
      \rounddown{(k+1)\,\eps}-b\le\rounddown{k\,\eps}
   \end{equation}
   holds for all $k\ge 1$, then $\eps\le b$. This latter
   assertion is
   obvious when $\eps$ is an integer. If $\eps=e+\delta$ with an
   integer $e$ and $0<\delta<1$, then there is an integer
   $k\ge 1$ such that
   $k\delta<1$
   and
   $(k+1)\delta\ge 1$.
   We then have
   $$
      \rounddown{(k+1)\eps}-b=(k+1)e+1-b
      \quad\mbox{and}\quad
      \rounddown{k\eps}=ke
   $$
   so that \eqnref{round-ineq} implies $\eps\le b$, as claimed.
\end{proofof}

\begin{remark}\rm
   The upper bound of $n+1$ in Corollary~\ref{cor-fano}
   follows also from the fact that for
   every point $x$
   on a smooth Fano variety $X$ of dimension $n$ there is a rational
   curve $C$ passing through $x$ and such that $-K_X\cdot C\leq n+1$.
   This is a deep fact proved by Mori and Koll\'ar (see \cite[Theorem
   V.1.6]{KolRC}). By contrast, our argument is fairly
   elementary.
\end{remark}

   As a further consequence of
   Proposition \ref{prop-fano}, we obtain the following characterization of
   $\P^n$ via Seshadri constants. Its statement will be strengthened
   considerably in Theorem~\ref{better_than_pn}.

\begin{corollary}\label{cor-Pn}
   Let $X$ be a smooth Fano variety such that
   $$
      \eps(-K_X,x)=n+1
   $$
   for all $x\in X$.
   Then
   $X\isom\P^n$.
\end{corollary}

\begin{proof}
   According to the previous result it is enough to show that
   the condition $\eps(-K_X,x)=\sqrt[n]{(-K_X)^n}=n+1$ implies
   that $X=\P^n$. This can be seen as follows. From the assumption
   and Proposition~\ref{prop-fano} we get
   $$
      s(k(-K_X),x)=k\,\eps(-K_X,x)=k(n+1)
   $$
   Putting $s=s(k(-K_X),x)$
   a result of Beltrametti and Sommese \cite[Theorem 3.1]{BS}
   implies then that
   \begin{equation}\label{BS-inequality}
      (k(-K_X))^n\ge s^n+s^{n-1}
   \end{equation}
   unless $X\isom\P^n$.
   Bearing in mind that
   $$
      (k(-K_X))^n=k^n(n+1)^n
   $$
   we get a contradiction with \eqnref{BS-inequality}, unless
   $X\isom\P^n$.
   (Note that the cited result in \cite{BS} assumes the line
   bundle in question to be \emph{$s$-jet ample}, whereas in our
   situation we only know that the bundle generates $s$-jets at
   all points. The proof in \cite{BS} works, however, under this
   weaker assumption.)
\end{proof}

\begin{remark}\rm
   For the proof of Corollary~\ref{cor-Pn}
   one could also invoke the following
   characterization of the projective space conjectured in
   \cite[Conjecture V.1.7]{KolRC} and proved by Kebekus in
   \cite{Keb02}:
   \begin{quote}\it
      A smooth projective variety $X$ of dimension $n$ is
      isomorphic to the projective space $\P^n$ if and only if
      it is Fano and $-K_X\cdot C\geq n+1$ for every rational curve
      $C\subset X$.
   \end{quote}
   The assumptions of the previous corollary imply that $-K_X$
   generates $n+1$ jets at every point $x\in X$ so that the
   inequality $-K_X\cdot C\geq n+1$ is fulfilled for an arbitrary curve $C$
   on $X$.
\end{remark}

   Now we show a considerably stronger version of
   Corollary~\ref{cor-Pn}. Of course this result implies
   Theorem~\ref{char-pn} stated in the introduction.

\begin{theorem}\label{better_than_pn}
   Let $X$ be a Fano variety such that
   $X\not\isom\P^n$,
   then the inequality
   $$
      \eps(-K_X,x)\le n
   $$
   holds for \textit{all} points $x\in X$.
\end{theorem}

   In view of
   Corollary~\ref{cor-fano} we need to show the
   following statement:
   \begin{quote}
      If $(-K_X)^n=(n+1)^n$ and if $\eps(-K_X,x)=n+1$ for some point
      $x\in X$, then $X\isom\P^n$.
   \end{quote}

   In the surface case
   it is easy to verify
   the above statement
   because we know exactly all Fano surfaces.

\begin{example}[Del Pezzo surfaces]\rm
   Let $X$ be a Del Pezzo surface, i.e., a smooth Fano variety of
   dimension two. Then
   \begin{equation}\label{eqn-Del-Pezzo}
      \eps(-K_X,x)\le 2
   \end{equation}
   for all $x\in X$, unless $X\isom\P^2$.
   In fact, $X$ is either
   \begin{itemize}
   \item[(i)]
      $\P^1\times\P^1$, or
   \item[(ii)]
      the blow-up of $\P^2$ in $r$ points, with
      $0\le r\le 8$.
   \end{itemize}
   In Case~(i) we have $K_X=\O(-2,-2)$ so that
   $\eps(-K_X,x)=2$ for all $x\in X$. In Case~(ii) we have
   $(-K_X)^2=9-r$, and therefore the number $\sqrt{(-K_X)^2}$
   is an integer only if $r=0$, $r=5$, or $r=8$. If $r=0$ then
   $X=\P^2$, and if $r=5$ or $r=8$ then $\sqrt{(-K_X)^2}=2$ or
   $\sqrt{(-K_X)^2}=1$ respectively.
   So we get the inequality \eqnref{eqn-Del-Pezzo} by applying
   Corollary~\ref{cor-fano}.
   Broustet \cite{Bro06} has recently determined the precise values of
   $\eps(-K_X,x)$ for
   $0\le r\le 8$.
\end{example}

   Now we give a proof of Theorem~\ref{better_than_pn} valid in arbitrary
   dimension. This proof is motivated by the methods of \cite{DRS00}.

\proofof{Theorem~\ref{better_than_pn}}
   As in the proof of \ref{cor-Pn} we get
   $$
      s((-K_X),x)=\eps(-K_X,x)=n+1
   $$
   in the fixed point $x$.
   Looking at the exact sequences defining bundles of $k$ jets of $-K_X$
   $$
     0\to \Sym^{k}\Omega_X\otimes(-K_X)\to J_k(-K_X)\to
     J_{k-1}(-K_X)\to 0
   $$
   for $k=1,\dots,n+1$ and computing inductively we obtain
   $$
     c_1(J_{n+1}(-K_X))=\O_X.
   $$
   Since by assumption the vector bundle $J_{n+1}(-K_X)$ is globally generated at
   the point $x$ and its determinant is trivial, it follows that the
   bundle itself is trivial. (This is because the determinant of global
   sections generating at $x$ does not vanish anywhere.)

   Now, let $f:\P^1\to X$ be a rational curve on $X$ (i.e., the
   map $f$ is non-constant). Let
   $$
     f^*(T_X)=\bigoplus_{i=1}^n\O(a_i) \mbox{ and } f^*(-K_X)=\O(b)\,.
   $$
   Note that $b>0$ since $-K_X$ is ample. Dualizing the defining
   exact sequence for $(n+1)$-st jets we have
   $$
     0\to J_n(-K_X)^*\to J_{n+1}(-K_X)^*\to
     \Sym^{n+1}T_X\otimes K_X\to 0\,.
   $$
   The bundle in the middle is trivial, it is in particular globally
   generated, hence the same is true for its quotient on the right.
   We write
   $$
      f^*(\Sym^{n+1}T_X)=\left(\bigoplus_{i=1}^n \O((n+1)a_i)\right)\oplus P\,,
   $$
   where $P$ abbreviates the remaining summands.
   Thus
   $$
      f^*(\Sym^{n+1}T_X\otimes K_X)=\left(\bigoplus_{i=1}^n \O((n+1)a_i-b)\right)\oplus
      P(-b)\,.
   $$
   It follows that
   $$
     (n+1)a_i-b\geq 0\,,
   $$
   which in view of $b>0$ implies $a_i>0$ for all $i=1,\dots,n$
   and we conclude by the Mori characterization of projective space
   \cite[Theorem V.3.2]{KolRC}.
\endproof

\section{Varieties with trivial canonical bundle}

   We consider now varieties whose canonical bundle is trivial.
   A straightforward modification of the proof of
   Proposition~\ref{prop-fano} yields the following statement:

\begin{proposition}\label{prop-K-zero}
   Let $X$ be a smooth projective variety of dimension $n$
   such that
   $K_X=\O_X$, let $L$ be an ample line bundle on $X$ and $x\in
   X$ a point.

   If $\eps(L,x)<\sqrt[n]{L^n}$ or if $\sqrt[n]{L^n}$ is not an
   integer, then one has for every integer $k\ge 1$ the
   inequalities
   $$
      \rounddown{k\,\eps(L,x)}-n\le
      s(kL,x)\le\rounddown{k\,\eps(L,x)} \ .
   $$
   In the alternative case \inparen{where
   $\eps(L,x)=\sqrt[n]{L^n}$ is an
   integer} one has
   $$
      k\sqrt[n]{L^n}-(n+1)\le s(kL,x)\le k\sqrt[n]{L^n} \,.
   $$
\end{proposition}

   So there are only $n+1$ potential values of $s(kL,x)$ in the
   first case, and $n+2$ potential values in the second case.
   This means that there is surprisingly little room for the
   numbers $s(kL,x)$.

\begin{example}\rm
   Consider a smooth quartic surface $X\subset\P^3$ containing a
   line $\ell$. Then for $L=\O_X(1)$ and $x\in\ell$ one has
   $\eps(L,x)=1$, and
   $$
      s(kL,x)=k=k\,\eps(L,x) \ .
   $$
   So in this case
   $s(kL,x)$ has its maximal possible value for all $k\ge 1$.
\end{example}

   The following proposition gives interesting
   constraints
   on maximal (in the sense of Proposition~\ref{prop-K-zero})
   values of the numbers $s(kL,x)$.

\begin{proposition}\label{cones}
   Let $X$ be a smooth projective variety, let $L$ be an ample
   line bundle on $X$, and let $x\in X$ be any point. If
   $\eps(L,x)$ is computed by a curve, i.e., if there is a curve
   $C\subset X$ such that
   $$
      \eps(L,x)=\frac{L\cdot C}{\mult_x(C)} \,,
   $$
   then
   \begin{itemize}
   \item[\rm(i)]
      $s(kL,x)<k\,\eps(L,x)$ for all $k\ge 1$, or
   \item[\rm(ii)]
      $C$ is smooth at $x$.
   \end{itemize}
   If $\dim(X)=2$, then in Case~(ii) one has $C^2\le 1$.
   If in addition $\eps(L,x)<\sqrt{L^2}$, then $C^2\le 0$.
\end{proposition}

   In other words, if one has $s(kL,x)=k\,\eps(L,x)$ for
   \textit{some}
   $k$, then the Seshadri constant $\eps(L,x)$ cannot be computed
   by a singular curve.

   Let us point out the following
   sample application of Proposition~\ref{cones}.

\begin{corollary}
   Let $X$ be an abelian surface of Picard number~one. Then for
   any ample line bundle $L$ on $X$, any $x\in X$,
   and any integer $k\ge 1$, one
   has
   $$
      s(kL, x)<k\,\eps(L,x)\,.
   $$
\end{corollary}

   This follows from the fact that in the situation of
   the corollary one knows from \cite[Sect.~6]{Bau99}
   that $\eps(L,x)$ is
   computed by a singular curve.

\proofof{Proposition~\ref{cones}}
   Let $k\ge 1$, and write $s=s(kL,x)$, $\eps=\eps(L,x)$, and
   $m=\mult_x(C)$. As $kL$ generates $s$-jets at $x$, there is a
   divisor $D\in|kL|$ with $\mult_x(D)=s$ and with prescribed
   tangent cone at the point $x$. So if $C$ is not smooth at $x$, then
   we can arrange that the projective tangent cones
   $$
      \P TC_x(D) \qquad\mbox{and}\qquad \P TC_x(C)
   $$
   intersect, while still $C\not\subset D$.
   Then by the intersection inequality
   \cite[Corollary~12.4]{Ful84}
   we have
   $$
      D\cdot C\ge s\cdot m+1 \,,
   $$
   and hence
   $$
      \eps=\frac{L\cdot C}m=\frac 1k\frac{D\cdot C}m\
      \ge\frac sk+\frac1{km} \,,
   $$
   which implies
   $$
      s\le\eps k-\frac 1m<\eps k \,.
   $$
   This proves the first assertion.

   Suppose now that $\dim(X)=2$, and that we are in Case~(ii).
   Then by the index theorem
   $$
      L^2C^2\le (L\cdot C)^2=\eps^2\le L^2 \,,
   $$
   and hence $C^2\le 1$. If $\eps<\sqrt{L^2}$, then the last
   inequality is strict, and we get $C^2\le 0$.
\endproof

   While the upper bound in Proposition~\ref{prop-K-zero} holds without any
   assumptions, the following example shows that one cannot expect
   a lower bound valid for all $k\geq 1$ without additional assumptions on
   the underlying variety.

\begin{example}\rm
   Let $k_0$ be a positive integer. We show that there exists a
   smooth projective
   surface $X$ and an ample line bundle $L$ on $X$ such that
   for all $k=1,\dots,k_0$ we have $s(kL,x)=-1$ for every point $x\in
   X$, whereas $\eps(L,x)=1$.

   To this end let $C$ be a curve of genus
   $g>k_0$ and let $D$ be a general divisor of degree $1$ on $C$.
   Then $h^0(C,D)=h^0(C,2D)=\dots=h^0(C,k_0D)=0$.
   Now let $X=C\times C$ be the product with projections $\pi_1$ and
   $\pi_2$. We set $L=\pi_1^*(D)\otimes \pi_2^*(D)$. The line bundle
   $L$ is ample and $\eps(L,x)=1$. By the K\"unneth formula we see
   that $h^0(X,L)=\cdots=h^0(X,k_0L)=0$ and hence $s(kL,x)=-1$
   for $k\le k_0$, as claimed.
\end{example}

   The next example shows that even in the surface case
   one cannot expect lower bounds on
   $s(kL,x)$ as in Proposition~\ref{prop-K-zero}
   without assumptions on
   $K_X$: The numbers $s(kL,x)$ may in fact
   be smaller than the lower
   bound of Proposition~\ref{prop-K-zero} for \textit{all} values
   of $k$.

\begin{example}\rm
   Let $X$ be a smooth surface of degree 9 in $\P^3$ with
   $\rho(X)=1$. For the line bundle $L=\O_X(1)$ one has
   $\eps(L,x)=\rounddown{\sqrt{L^2}}=3$ for very general $x\in
   X$,
   by Steffens result \cite{Ste98}.
   We assert that
   $$
      s(kL,x)<k\eps(L,x)-4 \,.
   $$
   for $k\gg 0$.
   In fact, if a line bundle
   generates $s$-jets at some point, then
   it must have
   at least $s+2\choose 2$ independent global
   sections. But by Riemann-Roch we have
   $$
      h^0(kL)=\chi(\O_X)+\frac 92 k(k-5) \,,
   $$
   and this is for $k\gg 0$
   not enough in order to generate jets of order
   $3k-1=k\,\eps(L,x)-4$.
   Note that the same argument works for any smooth surface
   $X\subset\P^3$ with $\rho(X)=1$, whose degree is a square
   number $\ge 9$.
\end{example}

   We now consider concrete
   applications of Proposition~\ref{prop-K-zero}.

\begin{example}[Jets of theta functions]\label{example-ppas}\rm
   Consider an irreducible principally polarized abelian surface
   $(X,\Theta)$. One knows
   that $\eps(\Theta,x)=\frac 43$ for
   every point $x\in X$
   (see \cite{Ste98} and \cite[Sect.~6]{Bau99}), hence
   Proposition~\ref{prop-K-zero} tells us that
   $$
      \rounddown{\tfrac 43 k}-2\le s(k\Theta,x)\le
      \rounddown{\tfrac 43 k}
   $$
   for every $k\ge 1$.
   So there are for each $k$ only three possible values that
   $s(k\Theta,x)$ can have. Here the possibilities
   are tabulated for small
   values of $k$:
   $$
   \renewcommand\arraystretch{1.2}
   \begin{array}{c|rrrrrrrrrrrrrrr}\hline
   k       & 1  & 2 & 3 & 4 & 5 & 6 & 7 & 8  & 9  & 10 \\ \hline
   s(kL,x) & 1  & 2 & 4 & 5 & 6 & 8 & 9 & 10 & 12 & 13 \\
           & 0  & 1 & 3 & 4 & 5 & 7 & 8 & 9  & 11 & 12 \\
           & -1 & 0 & 2 & 3 & 4 & 6 & 7 & 8  & 10 & 11 \\ \hline
   \end{array}
   $$
   Of course it is very hard to determine what the exact values are
   -- especially as, contrary to the value of
   $\eps(\Theta,x)$, they might and do depend on the point $x$.
   Specifically, as
   $|\Theta|=\{\Theta\}$, one has
   $s(\Theta,x)=0$ if $x\notin\Theta$, and $s(\Theta,x)=-1$ if
   $x\in\Theta$.
   As for $s(2\Theta,x)$, one knows that the linear series
   $|2\Theta|$ defines a map $X\to\P^3$ of degree 2 onto the
   Kummer surface, hence
   $s(2\Theta,x)\ge 1$ for
   generic $x$, but
   $s(2\Theta,x)=0$ for points mapped onto double
   points of the Kummer quartic in $\P^3$.

   Counting sections carefully
   we can in fact rule out several values in the above table. Below
   we present the remaining possibilities. Values marked in bold are
   actually taken on.
   $$
   \renewcommand\arraystretch{1.2}
   \begin{array}{c|rrrrrrrrrrrrrrr}\hline
   k       & 1           & 2          & 3          & 4 & 5 & 6 & 7 & 8  & 9  & 10 \\ \hline
   s(kL,x)        & \textbf{0}  & \textbf{1} &            & 4 & 5 & 7 & 8 & 9  & 11 & 12 \\
           & \textbf{-1} & \textbf{0} & \textbf{2} & 3 & 4 & 6 & 7 & 8  & 10 & 11 \\ \hline
   \end{array}
   $$

   There are (at least) two things that would be
   interesting to know in this context:
   \begin{itemize}
   \item
      Is $s(k\Theta,x)$ independent of $x$ when $k\gg 0$ ?
   \item
      Is the sequence of numbers $s(k\Theta,x)$ (for large $k$) the same for every
      principally polarized abelian surface $(X,\Theta)$ or does
      it depend on the moduli?
   \end{itemize}
\end{example}

   In the following two examples, yet more precise statement about
   the numbers $s(kL,x)$ can be obtained.

\begin{example}[Abelian surfaces of type $(1,2)$]\rm
   It can happen for certain line bundles that $s(kL,x)$ is not
   only below $k\,\eps(L,x)$, as predicted by
   Proposition~\ref{cones}, but even below $k\,\eps(L,x)-1$ for
   \textit{all} values of $k$. Consider for instance an abelian
   surface $X$ of Picard number 1 carrying a polarization $L$ of
   type $(1,2)$. Then $L^2=4$ and, by
   \cite[Sect.~6]{Bau99}, $\eps(L)=2$.
   A line bundle that generates $s$-jets at some point must have
   at least $s+2\choose 2$ independent global
   sections. As we have
   by Riemann-Roch
   $h^0(kL)=2k^2$, this implies that $kL$ cannot generate
   $2k$-jets or $(2k-1)$-jets at any point. So
   $$
      s(kL,x)=2k-3 \quad\mbox{or}\quad
      s(kL,x)=2k-2
   $$
   for every $x\in X$.
\end{example}

\begin{example}[Quartic surfaces]\rm
   Let $X$ be a smooth quartic surface in $\P^3$ and
   $L=\O_X(1)$. For general $x\in X$ one has $\eps(L,x)=2$ by \cite[Theorem]{Bau97}.
   As $h^0(kL)=2+2k^2$, the
   dimension argument of the previous example
   gives the same conclusion as there:
   $$
      s(kL,x)=2k-3 \quad\mbox{or}\quad
      s(kL,x)=2k-2 \,.
   $$
\end{example}

   The two previous examples are special cases of:

\begin{proposition}
   Let $X$ be a smooth projective surface with $K_X=\O_X$, and
   let $L$ be an ample line bundle on $X$.
   Suppose that $\eps(L,x)=\sqrt{L^2}$ at a fixed point $x\in X$.

   If $\sqrt{L^2}$ is an integer, then
   we have for $k\gg 0$
   $$
      s(kL,x)=k\sqrt{L^2}-3 \quad\mbox{or}\quad
      s(kL,x)=k\sqrt{L^2}-2 \,.
   $$
   and in the alternative case we have
   $$
      s(kL,x)=\rounddown{k\sqrt{L^2}}-2 \,,
   $$
   for infinitely many values of $k$.
\end{proposition}

\begin{proof}
   For the linear series
   $|kL|$ to generate $s$-jets at $x$, the line bundle $kL$ needs
   to have
   at least $s+2\choose 2$ independent global
   sections.
   By the Riemann-Roch theorem we have for $k\gg 0$
   $$
      h^0(kL)=\chi(\O_X)+\frac12{L^2} k^2.
   $$
   If $\eps(L,x)=\sqrt[n]{L^n}$ is an integer, then
   one finds with elementary calculations that
   the inequality
   $h^0(kL)\ge {s+2\choose 2}$ cannot hold for
   $k\gg 0$ when
   $s\geq \rounddown{k\,\eps(L,x)}-1$.
   The assertion follows then from Proposition~\ref{prop-K-zero}.
   In the remaining case one has to work with the rounddown and
   therefore the assertion is somewhat weaker.
\end{proof}

\paragraph*{Acknowledgement.}
   We would like to thank Jaros\l{}aw Wi\'sniewski for helpful
   discussions and for bringing \cite{DRS00} to our attention.
   This paper was written while the second author was visiting
   professor at the University of Duisburg-Essen. It is a pleasure
   to thank H\'el\`ene Esnault and Eckart Viehweg for their invitation
   and hospitality. The second author was partially supported by
   MNiSW grant N~N201 388834.

\bigskip
\small
   Tho\-mas Bau\-er,
   Fach\-be\-reich Ma\-the\-ma\-tik und In\-for\-ma\-tik,
   Philipps-Uni\-ver\-si\-t\"at Mar\-burg,
   Hans-Meer\-wein-Stra{\ss}e,
   D-35032~Mar\-burg, Germany.

\nopagebreak
   \textit{E-mail address:} \texttt{tbauer@mathematik.uni-marburg.de}

\bigskip
   Tomasz Szemberg,
   Instytut Matematyki UP,
   PL-30-084 Krak\'ow, Poland

\nopagebreak
   \textit{E-mail address:} \texttt{szemberg@ap.krakow.pl}

\medskip
   \textit{Current address:}
   Instytut Matematyczny PAN, ul. \'Sniadeckich 8, PL-00-956 Warszawa, Poland


\begin{thebibliography}{99}\footnotesize\itemsep=0cm\parskip=0cm

\bibitem{Bau97}
   Bauer, Th.:
   Seshadri constants of quartic surfaces.
   Math. Ann.  309  (1997), 475--481

\bibitem{Bau99}
   Bauer, Th.:
   Seshadri constants on algebraic surfaces.
   Math. Ann. 313 (1999), 547--583

\bibitem{BS}
   Beltrametti, M.C., Sommese, A.J.:
   On $k$-jet ampleness.
   Complex analysis and geometry (Ancona, Vincenzo, eds.), New York, Plenum Press, 1993, pp. 355--376

\bibitem{Bro06}
   Broustet, A.:
   Constantes de Seshadri du diviseur anticanonique des surfaces de del Pezzo.
   L'Enseignement Math. 52 (2006), 231--238

\bibitem{Dem92}
   Demailly, J.-P.:
   Singular Hermitian metrics on positive line bundles.
   Complex algebraic varieties (Bayreuth, 1990), Lect. Notes Math. 1507,
   Springer-Verlag, 1992, pp. 87--104

\bibitem{DRS00}
   Di Rocco, S., Sommese, A.J.:
   Line bundles for which a projectivized jet bundle is a product.
   Proc. Amer. math. Soc. 129 (2000) 1659--1663

\bibitem{Ful84}
   Fulton, W.:
   Intersection theory.
   Ergeb. Math. Grenzgeb. (3) 2, Springer-Verlag, 1984

\bibitem{Keb02}
   Kebekus, S.:
   Characterizing the projective space after Cho, Miyaoka and
   Shepherd-Barron. In Complex geometry (Gottingen, 2000), 147--155,
   Springer, Berlin, 2002.

\bibitem{KolRC}
   Koll\'ar, J.:
   Rational curves on algebraic varieties.
   Ergebnisse der Mathematik und ihrer Grenzgebiete. 32.
   Springer-Verlag, Berlin, 1996

\bibitem{PAG}
   Lazarsfeld, R.:
   Positivity in Algebraic Geometry~I.
   Springer-Verlag, 2004.

\bibitem{Ste98}
   Steffens, A.:
   Remarks on Seshadri constants.
   Math. Z. 227 (1998), 505--510

\end{thebibliography}
\end{document}